\documentclass{amsart}
\usepackage{amssymb,amsfonts}
\usepackage[all,arc]{xy}
\usepackage{enumerate}
\usepackage{mathrsfs}
\usepackage{microtype}
\usepackage{multicol}
\usepackage{xcolor}
\usepackage{upgreek}
\usepackage{graphicx} 
\usepackage{cancel,shuffle}
\usepackage{stackrel}
\usepackage{enumerate}
\usepackage{graphicx}
\usepackage{fdsymbol}
\usepackage[all]{xy}
\usepackage{tikz-cd}
\usetikzlibrary{cd}
\usepackage{hyperref}

\newtheorem{theorem}{Theorem}[section]
\theoremstyle{definition}
\newtheorem{example}[theorem]{Example}

\newcommand{\RR}{\mathbb{R}}
\newcommand{\e}{\mathtt{e}}
\newcommand{\letd}{\mathtt{d}}
\newcommand{\wi}{\mathtt{i}}
\newcommand{\w}{\mathtt{w}}
\newcommand{\lv}{\mathtt{v}}
\newcommand{\la}{\mathtt{a}}
\newcommand{\lb}{\mathtt{b}}
\newcommand{\halfshuffle}{\mathbin{\succ}}

\makeatletter
\let\c@equation\c@theorem
\makeatother
\numberwithin{equation}{section}

\begin{document}

\title{Computing Path Signature Varieties in Macaulay2}
\author{Carlos Am\'endola}
\email{amendola@math.tu-berlin.de}

\author{Angelo El Saliby}
\email{angelo.elsaliby@studenti.units.it}

\author{Felix Lotter} 
\email{felix.lotter@mis.mpg.de}

\author{Oriol Reig Fité}
\email{oriol.reigfite@unitn.it}



\begin{abstract}
    The signature of a path is a non-commutative power series whose coefficients are given by certain iterated integrals over the path coordinates. This series almost uniquely characterizes the path up to translation and reparameterization. Taking only fixed degree parts of these series yields signature tensors. We introduce the \textsc{Macaulay2} package \texttt{PathSignatures} to simplify the study of these interesting objects for piecewise polynomial paths. It allows for the creation and manipulation of parametrized families of paths and provides methods for computing their signature tensors and their associated algebraic varieties.
\end{abstract}

\maketitle

\section{Introduction}

Given a path $X:[0,1]\to \mathbb R^d$ defined by $d$ piecewise differentiable functions $X_1, \dots, X_d :[0,1]\rightarrow \mathbb{R}^d$, its \textit{signature} is a linear form $\sigma(X): T\bigl((\mathbb{R}^d)^*\bigr)\rightarrow \mathbb{R}$ on the tensor algebra of the dual of 
            $\mathbb{R}^d$, whose image on a decomposable tensor $\alpha_1\otimes \dots\otimes \alpha_k$ is the iterated integral $$
            \alpha_1\otimes \dots\otimes \alpha_k\overset{\sigma}{\mapsto} \int_0^1\int_0^{t_k}\dots\int_0^{t_2} (\alpha_1 X)'(t_1)\dots (\alpha_k X)'(t_k) \mathrm{d} t_1\dots \mathrm{d}t_k.
            $$
The \textit{$k$-th level signature} $\sigma^{(k)}(X)$ is the restriction of $\sigma(X)$ to $\bigl((\mathbb R^d)^*\bigr)^{\otimes k}$. Since the dual of $\bigl((\mathbb R^d)^*\bigr)^{\otimes k}$ is canonically isomorphic to $(\mathbb{R}^d)^{\otimes k}$, the $k$-th level signature is also called the $k$-th \textit{signature tensor} of $X$.\par

While path signatures have long had prominence in stochastic analysis \cite{lyons2014rough}, they have been more recently approached from an algebraic statistics and nonlinear algebra perspective \cite{sullivant2023algebraic,michalek2021invitation}. A key point is that one can define algebraic varieties from signature tensors \cite{AFS18, galuppi2019rough, colmenarejo2020toric} and study the way signatures are transformed by polynomial maps \cite{colmenarejo2020}. As tensors, one can also naturally study questions such as their rank \cite{amendola2025, galuppi2024rank} and use this algebraic structure for path learning \cite{pfeffer2019learning}.

Even though computer systems like \textsc{Macaulay2} \cite{grayson2002macaulay2} and \textsc{Bertini} \cite{bates2013numerically} have been used to study path signatures \cite{AFS18}, there is still no package providing the algebraic structures and methods relevant in this context. We present a package in \textsc{Macaulay2} that aims to fill this gap and demonstrate how it simplifies the study of signature tensors with several key examples.

\section{Piecewise Polynomial Paths}
We focus on piecewise polynomial paths, that is, paths $X:[0,1]\to\mathbb R^d$ with the property that there are $0 = t_0 < \ldots < t_m =1$ such that the coordinate functions of the restrictions $X_i := X|_{[t_i,t_{i+1}]}$ are given by polynomials. After reparametrization and translation, each $X_i$ is a polynomial path $[0,1]\to\mathbb R^d$ starting at the origin. As the signature $\sigma(X)$ is invariant under both reparametrization and translation, we are only interested in paths up to these equivalence relations; thus, we can encode the piecewise polynomial path $X$ by its \textit{polynomial segments} $X_0,\dots,X_m$. The type \texttt{Path} records precisely this data. It includes as special cases the families of piecewise linear paths and polynomial paths studied in \cite{AFS18}. \par
A key feature of the package is that it also allows to consider \textit{parametrized families} of piecewise polynomial paths: the polynomial segments defining a path are allowed to have coefficients in any ring. In the following example, we illustrate this for a linear path.

\begin{example}[Linear paths]\label{ex:linpath}
    Piecewise linear paths are a classical way to approximate general paths. A linear path can be constructed with the method \texttt{linPath} by giving an increment.

\begin{verbatim}
    i1: loadPackage "PathSignatures";
    i2: QQ[x_1, x_2]; -- The increments can be in any ring
    i3: Z = linPath({2*x_1,3*x_2}) 
    o3 = Path in 2-dimensional space with 1 polynomial segment:

         {{2x t, 3x t}}
             1     2

    o3 : Path
\end{verbatim}
A piecewise linear path can be constructed by concatenation of linear paths (as described in Example \ref{ex:monpath}), or with the method \texttt{pwLinPath} by giving a matrix of increments.

\begin{verbatim}
    i4 : M = id_(QQ^3);
    i5 : X = pwLinPath(M)
    o5 = Path in 3-dimensional space with 3 polynomial segments:

         {{t, 0, 0}, {0, t, 0}, {0, 0, t}}

    o5 : Path
\end{verbatim}
\end{example}

\begin{example}\label{ex:monpath}
A general polynomial path can be created with the method \texttt{polyPath}, by giving a list containing its coordinate functions. Concatenating them with the symbol \texttt{**}, one obtains piecewise polynomial paths.
\begin{verbatim}
    i6 : R = QQ[t];
    i7 : Y = polyPath({t, t^2, t^3})
    o7 = Path in 3-dimensional space with 1 polynomial segment:

               2   3
         {{t, t , t }}
     
    o7 : Path
    i8 : X**Y
    o8 = Path in 3-dimensional space with 4 polynomial segments:

                                                 2   3
         {{t, 0, 0}, {0, t, 0}, {0, 0, t}, {t, t , t }}

    o8 : Path
\end{verbatim}
\end{example}

Note that according to our encoding of paths, concatenation is only formal, and no reparameterization is attempted. The path signature is computed from the polynomial pieces through Chen's Formula (\cite[Corollary~5.1]{AFS18}). 

\section{Free algebras of words, Lyndon Words and Shuffles}
Consider the $d$ letters $\mathtt{1},...,\letd$ and denote by $\mathbb R\langle \mathtt{1},...,\letd\rangle$ the free associative algebra over them. Then the algebra $T\bigl((\mathbb{R}^d)^*\bigr)$ may be identified with $\mathbb R\langle \mathtt{1},...,\letd\rangle$ via the isomorphism induced by $e_i^* \mapsto \mathtt{i}$. A \emph{word} is a concatenation of letters, i.e., a product of letters in the free associative algebra.\par 

The method $\texttt{wordAlgebra(d)}$ creates the free associative algebra over the $d$ letters  $\mathtt{Lt_1}, \dots, \mathtt{Lt_d}$. The resulting ring is of type \texttt{NCPolynomialRing}, which is provided by the \textsc{Macaulay2} package \texttt{NCAlgebra}. A linear combination of words is encoded by the type \texttt{NCRingElement}.\par
A useful notation for words, implemented in the method \texttt{wordFormat}, is the following: to each word $\mathtt{i_1}\dots \mathtt{i_k}$ we associate the array $[\mathtt{i_1}, \dots, \mathtt{i_k}]$ and then extend this convention by linearity to all elements of $\mathbb{R}\langle\mathtt{1},\dots\mathtt{d}\rangle$ by considering formal linear combinations of such arrays. This notation can also be used through the method \texttt{Array}~\textunderscore\,~\texttt{NCPolynomialRing } to define words in an algebra.
\begin{verbatim}
    i9  : A = wordAlgebra(2);
    i10 : w = [1,2]_A
    
    o10 = Lt Lt
            1  2
            
    o10 : A
\end{verbatim}

Taking the path $Z$ in Example \ref{ex:linpath}, its signature $\sigma(Z)$ is represented by the method \texttt{sig}. For instance it can be evaluated at a word:

\begin{verbatim}
    i11 : sig(Z,w)
    
    o11 = 3x x
            1 2

    o11 : QQ[x ..x ]
              1   2

\end{verbatim}
Or we can extract its $\mathtt{k}$-th level signature tensor:
\begin{verbatim}
    i12 : sig(Z, 2) // wordFormat
    
          9 2                                          2
    o12 = -x  [2, 2] + 3x x  [2, 1] + 3x x  [1, 2] + 2x  [1, 1]
          2 2            1 2            1 2            1
    
\end{verbatim}
With the notation established above the last output corresponds to the tensor $$\frac{9}{2} x_2^2 \cdot \mathtt{Lt_2}^2 + 3 x_1 x_2 \cdot \mathtt{Lt_2 Lt_1} +3 x_1 x_2 \cdot \mathtt{Lt_1 Lt_2} + 2 x_1^2 \cdot \mathtt{Lt_1}^2.$$
We denote the usual product of $\mathbb R \langle \mathtt{1},\ldots,\mathtt{d} \rangle$, which is just concatenation on words, by $\smallblackcircle$. We write $\e$ for the identity element of the multiplication, the \textit{empty word}.\par
There is another algebra structure on $\mathbb R \langle \mathtt{1},\ldots,\mathtt{d} \rangle$, known as the \emph{shuffle} product. It encodes a combinatorial relation in the signatures values. Intuitively, the shuffle product of two words is the sum over all words obtained by \emph{interleaving} them, with multiplicity. More precisely, if $\w,\w_1,\w_2$ are words and $\la,\lb$ are letters, it can be defined recursively as
\begin{eqnarray*}
\begin{array}{l}
\e \shuffle \w = \w \shuffle \e = \w, \\
(\w_1\smallblackcircle \la) \shuffle (\w_2\smallblackcircle \lb) = \bigl(\w_1\shuffle (\w_2\smallblackcircle \lb)\bigr)\smallblackcircle \la + \bigl((\w_1\smallblackcircle \la )\shuffle \w_2\bigr)\smallblackcircle\lb.
\end{array}
\end{eqnarray*}

Note that the shuffle product is commutative. For any path $X$, the signature $\sigma(X): \RR\langle \mathtt{1},...,\letd\rangle_{\shuffle}\to \RR$ is an algebra homomorphism (see \cite{Ree}).\par
Using \texttt{PathSignatures}, one can compute the shuffle product of two words 
with the method \texttt{NCRingElement\,\,**\,\,NCRingElement}.\par

\begin{example}
    Let us check that $\sigma(X)$ is a shuffle algebra homomorphism in an example:
    \begin{verbatim}
    i1: loadPackage "PathSignatures";
    i2 : A = wordAlgebra(3);
    i3 : R = QQ[t];
    i4 : X = polyPath({t,t^2,t^3});
    i5 : a = [1,2]_A;
    i6 : b = [2,3]_A;
    i7 : c = a ** b; c // wordFormat
    o8 = [2, 3, 1, 2] + [2, 1, 3, 2] + [2, 1, 2, 3] 
         + [1, 2, 3, 2] + 2 [1, 2, 2, 3]
    i9 : sig(X,a) * sig(X,b) == sig(X,c)
    o9 = true   
    \end{verbatim}
\end{example}

The shuffle product can also be defined as the symmetrization of the \emph{half shuffle} product, defined on $T^{\geq 1}(\RR^d)$, the space spanned by nonempty words, recursively as
 \begin{align*}
  \w\halfshuffle\wi &:= \w\smallblackcircle\wi,\\
  \w\halfshuffle (\lv\smallblackcircle\wi) &:= (\w\halfshuffle\lv+\lv\halfshuffle\w)\smallblackcircle\wi,
 \end{align*}
where $\w,\lv$ are words and $\wi$ is a letter. A direct calculation shows that $\w \shuffle \lv=\w\halfshuffle \lv + \lv\halfshuffle \w$. The method \texttt{NCRingElement~>>~NCRingElement} computes the half shuffle of two words.\par
The half shuffle plays an important role in translating between properties of paths and properties of signatures (see \cite{preiss2023algebraic}). For example, it is used to describe the transformations of the signature corresponding to polynomial transformations of paths, as we will explain in Section \ref{sec:transformations}.

\begin{example}
    Let us play with shuffles and half shuffles of two words:
    \begin{verbatim}
    i10 : R = wordAlgebra(3);
    i11 : w = [1]_R;
    i12 : v = [1,2,3]_R;
    i13 : s = w ** v; --shuffle product of w, v
    i14 : s // wordFormat
    o14 = [1, 2, 3, 1] + [1, 2, 1, 3] + 2 [1, 1, 2, 3]
    i15 : (v >> w) // wordFormat
    o15 = [1,2,3,1]
    i16 : (w >> v) + (v >> w) == w ** v
    o16 = true
    i17 : (w >> (v >> s)) == ((w ** v) >> s)
    o17 = true
    \end{verbatim}
    The last two lines verify an instance of the \textit{Zinbiel identity} of the non-associative half-shuffle algebra (\cite[Theorem~5]{colmenarejo2020}).
\end{example}

Since the algebra $\mathbb R \langle \mathtt{1}, \ldots, \mathtt{d} \rangle_\shuffle$ is commutative, it can be represented as a quotient of an (infinitely generated) polynomial ring. In fact, it turns out that it is a free commutative algebra over a suitable set of generators; one such choice are \textit{Lyndon words}.\par
 A word $l$ on the alphabet $\{\texttt{1},\dots, \texttt{d}\}$ is a {\em Lyndon word} if it is strictly smaller, in lexicographic order, than all of its rotations (or, equivalently, if it is smaller than all of its proper right subwords). The method \texttt{lyndonWords} computes all Lyndon words of at most a given length on a given number of letters:

\begin{verbatim}
    i18 : lyndonWords(3,2)
    o18 = {[1], [1, 2], [1, 3], [2], [2, 3], [3]}
    o18 : List
\end{verbatim}

Given an element of $\mathbb R \langle \mathtt{1}, \ldots, \mathtt{d} \rangle_\shuffle$,  the method \texttt{lyndonShuffle} computes its representation as a shuffle polynomial in Lyndon words. This is encoded as a \texttt{hashTable} mapping monomials (also given as a \texttt{hashTable}) to their coefficients.
\begin{example}
    We compute the representation of $\texttt{321}$ in $\mathbb R \langle \mathtt{1}, \mathtt{2}, \mathtt{3} \rangle$ in terms of Lyndon words.
\begin{verbatim}
    i19 : A = wordAlgebra(3);
    i20 : w = [3,2,1]_A;
    i21 : lyndonShuffle w
    o21 = HashTable{HashTable{[1, 2, 3] => 1} => 1}
                   HashTable{[1, 2] => 1} => -1
                             [3] => 1
                   HashTable{[1] => 1   } => -1
                             [2, 3] => 1
                   HashTable{[1] => 1} => 1
                             [2] => 1
                             [3] => 1
    o21 : HashTable
 \end{verbatim}
 That is, $$\texttt{321} = \texttt{123} - \texttt{12} \shuffle \texttt{3} - \texttt{1} \shuffle \texttt{23} + \texttt{1} \shuffle \texttt{2} \shuffle \texttt{3}.$$
\end{example}

\bigskip

\section{Universal varieties}

To any Lyndon word $l$ we can associate an iterated Lie bracketing $b(l) \in \mathbb R\langle \mathtt{1},\ldots,\mathtt{d}\rangle$, defined iteratively as follows. For a letter $\mathtt{i}\in \{\mathtt{1},\dots, \mathtt{d}\}$
we simply define $$ b(\mathtt{i}) := \mathtt{i}.$$
For a Lyndon word $l$ of length greater than 1 we define $$
b(l) := [b(l_1), b(l_2)],$$
where $l_1, l_2$ are non-empty words such that their concatenation $l_1 l_2$ is $l$, and $l_2$ is the longest Lyndon word appearing as a proper right factor of $l$. 

\bigskip

One can show that the bracketings $b(l)$ of Lyndon words $l$ of length at most $k$ form a basis for the Lie algebra $\text{Lie}^k(\RR^d)$ (see \cite[Theorem 4.9, Theorem 5.1]{reutenauerfree93}), justifying the name  \texttt{lieBasis} for the method computing the bracketing of a Lyndon word. 

\bigskip

The corresponding Lie group $\exp(\text{Lie}^k(\RR^d))$ is an affine variety, and the Chen-Chow theorem states that its projection to the level $k$ component describes the set of all possible $k$-th level signature tensors of smooth paths $X:[0,1]\rightarrow\RR^d$ (see \cite{Friz_Victoir_2010}). Inspired by this fact, the complexification and projectivization of this projection onto level $k$ is known as the \emph{universal variety} $\mathcal{U}_{d,k}$ (\cite[Section 4]{AFS18}). 

\begin{example}
    Let $d=2$ and $k=3$. The following code computes the parametrization $p_k \circ \exp: \text{Lie}^k(\RR^d) \to \mathcal{U}_{d,k}$ of the universal variety, where $p_k$ denotes the projection to the level $k$ component. First, we compute a basis of the Lie algebra as described above:
\begin{verbatim}
    i1 : loadPackage "PathSignatures";
    i2 : words = toList apply((3:1)..(3:2), i-> new Array from i);
    i3 : lwords = lyndonWords(2,3)
    o3 = {[1], [1, 1, 2], [1, 2], [1, 2, 2], [2]}
    o3 : List
    i4 : Q = QQ new Array from (apply(lwords,i->y_i) |
             {Degrees => apply(lwords, i->length(i))});
    i5 : A2 = wordAlgebra(2, CoefficientRing => Q);
    i6 : lbasis = apply(lwords, i-> lieBasis(i,A2));
\end{verbatim}
The variables in \texttt{Q} will represent the coordinates of an element of $\text{Lie}^k(\RR^d)$ with respect to this basis; in other words, we write an arbitrary element of $\text{Lie}^k(\RR^d)$ as
\begin{verbatim}
    i7 : lT = sum(0..length(lbasis)-1, i-> Q_i * lbasis_i);
\end{verbatim}
Its image under $p_k \circ \exp$ is computed by the method \texttt{tensorExp}:
\begin{verbatim}
    i8 : gT = tensorExp(lT,3);
\end{verbatim}
This is a non-commutative polynomial with coefficients in $\texttt{Q}$, or equivalently a tensor whose entries are polynomials in $\texttt{Q}$. We obtain the ring map corresponding to this polynomial parametrization as follows:
\begin{verbatim}
    i9 : m = tensorParametrization(gT);
    o9 : RingMap Q <-- QQ[...]
\end{verbatim}
We can then compute the prime ideal of the universal variety $\mathcal{U}_{2,3}$ by solving the corresponding implicitization problem:
\begin{verbatim}
    i10 : I = ker m;
    o10 : Ideal of QQ[...]
    i11 : dim I
    o11 = 5
    i12 : degree I
    o12 = 4
    i13 : betti mingens I
    
                 0 1
    o13 = total: 1 6
              0: 1 .
              1: . 6
\end{verbatim}    
Thus, $\mathcal{U}_{2,3}$ has affine dimension 5, degree 4 and is generated by 6 quadrics. This is consistent with \cite[Table 2]{AFS18}.
\end{example}

\section{Path signatures under linear and polynomial transformations}\label{sec:transformations}

From the definition of the signature it is immediate that it is equivariant under linear transformations in the following sense:
if $X: [0,1] \to \mathbb R^d$ is a path and $A$ is an $e\times d$ matrix, then
$$\sigma(A \circ X) = \sigma(X) \circ T(A^\top)$$
where $T(A^\top)$ is the map induced by $A^\top$ on the tensor algebra. In other words, using the canonical pairing $\bigl(\mathbb R\langle\texttt{1},\ldots,\texttt{d} \rangle \bigr)^* \times \bigl(\mathbb R\langle\texttt{1},\ldots,\texttt{d} \rangle \bigr)$, we have
$$\bigl\langle \sigma(A\circ X)\,,\, \mathtt{w} \bigr\rangle = \bigl\langle \sigma(X)\,,\, A^\top .w\bigr\rangle,$$
where the action of the matrix on the tensor, denoted by $A^T.w$, is diagonal.
\begin{example}
    Consider the following polynomial path \begin{verbatim}
    i1 : loadPackage "PathSignatures";
    i2 : R = QQ[t]; 
    i3 : Z = polyPath({t+2*t^2+3*t^3, 4*t +5*t^2+6*t^3})
    o3 = Path in 2-dimensional space with 1 polynomial segment:
    
             3     2        3     2
         {{3t  + 2t  + t, 6t  + 5t  + 4t}}
         
    o3 : Path
    \end{verbatim}
    Then $\texttt{Z}$ can be obtained as the image of the path $\texttt{Y}$ of Example \ref{ex:monpath} through the linear transformation $A : \mathbb{R}^3\rightarrow \mathbb{R}^2$ determined by the matrix \begin{verbatim}
    i4 : A = matrix({{1,2,3},{4,5,6}})
    
    o4 = | 1 2 3 |
         | 4 5 6 |
         
                   2       3
    o4 : Matrix ZZ  <-- ZZ
    \end{verbatim}
    Through the method \texttt{Matrix\,*\,NCRingElement} we can compute the diagonal action of a matrix on a tensor. As an example we verify the formula above.\begin{verbatim}
    i5 : Y = polyPath({t, t^2, t^3});
    i6 : T = wordAlgebra(2);
    i7 : sig(Z, Lt_1*Lt_2)  
    
         427
    o7 = ---
          10
           
    o7 : QQ
    
    i8 : sig(Y,(transpose A)*( Lt_1*Lt_2))
    
         427
    o8 = ---
          10
           
    o8 : QQ
    \end{verbatim} 
    We will generalize this example in Section \ref{section:core-signature-tensors-and-varieties}.
\end{example}

In \cite{colmenarejo2020}, the equivariance relation was generalized to polynomial transformations. Let $X$ be a path in $\RR^d$ with $X(0)=0$ and let $p:\RR^d\to \RR^m$ be a polynomial map such that $p(0)=0$. Then there is an algebra homomorphism $M_p:(\RR^m, \shuffle)\to (\RR^d, \shuffle)$ such that, for every word $\w\in T(\RR^m)$, we have
\begin{equation}\label{Eq: adjoint}
    \bigl\langle \sigma(p \circ X)\,,\, \w\bigr\rangle =\bigl\langle \sigma(X)\,,\, M_p(\w)\bigr\rangle.
\end{equation}

Moreover, consider the algebra homomorphism 
\begin{eqnarray*}
\begin{array}{rrcl}
\varphi_d: & \mathbb{R}[x_1,x_2,\dots,x_d] & \longrightarrow & \left(T(\mathbb{R}^d),\shuffle \right) \\
& x_i & \longmapsto & \wi \\
& x_{i_1}\cdots x_{i_l} & \longmapsto & \wi_1\shuffle \dots \shuffle \wi_l
\end{array}
\end{eqnarray*}
Then the restriction of $M_p$ to $T^{\geq 1}(\RR^d)$ is the unique half-shuffle homomorphism such that
\[M_p(\wi)=\varphi_d(p_i), \,\forall i\in\{1, \dots, d\}.\]
See \cite{colmenarejo2020} for more on these homomorphisms and for the relevant proofs.

\begin{example}
    Let us test the adjoint relation of \ref{Eq: adjoint}. The method \texttt{adjointWord} computes the evaluation of the corresponding map $M_p$. 
\begin{verbatim}
   i9 : S = QQ[x,y];
   i10 : p = {x^2,y^3,x-y};  
   i11 : R = QQ[t];
   i12 : X = polyPath({t,t^2});
   i13 : PP = apply(p, q -> sub(q, {x=>t, y=>t^2}));
   i14 : Y = polyPath(PP);
   i15 : wA3 = wordAlgebra(3);
   i16 : sig(Y,Lt_1)==sig(X, adjointWord(Lt_1, wA3,p))

   o16 = true
        
\end{verbatim}

    % f = [1,2]_R; -- [i_1,...,i_k]_R defines a word
    %         wordFormat f -- display the polynomial in word notation
    %         f = ([1,2]_R ** [1,2]_R); -- compute the shuffle product
    %         f // wordFormat --display the result in word notation

\end{example}

\section{Core Signature Tensors and Path Signature Varieties}\label{section:core-signature-tensors-and-varieties}
Given a path $X$ in $\mathbb{R}^d$, the adjoint relation in (\ref{Eq: adjoint}), when the polynomial map $p:\RR^d\to \RR^m$ is linear, becomes 
\[\sigma(p\circ X)=\sigma(X)\circ T(p^*), \]
where $T(p^*)$ is the induced linear map on the tensor algebra. More concretely, we have $\sigma^{(k)}(p\circ X)=A .\sigma^{(k)}(X)$ for every $k\geq 0$, where $A$ is the matrix corresponding to $p$, and it acts diagonally in the canonical way on $\sigma^{(k)}(X)$. In other words, the signature map is \emph{equivariant} under linear transformations.

\medskip

This leads us to consider the following setting: we have a set $\mathcal{S}$ of paths in $\RR^d$, and assume that there exists a path $D:[0,1]\to \RR^d$ such that 
\[\mathcal{S}=\{A\circ D \;| \; A:\RR^m\to \RR^d \; \text{linear}\}.\]
Then, to determine the signature of any path in $\mathcal{S}$,  we only need to compute the signature of $D$ and use the equivariance relation. We call the signature of $D$ at level $k$ the \textit{$k$-th core tensor} of the family $\mathcal{S}$.

\medskip

There are two instances of this situation that arise typically. The first is to take $\mathcal{S}$ as the set of piecewise linear paths made of $d$ linear steps. We  define the {\em canonical axis path} in $\mathbb{R}^d$ 
to be the path from $(0, \dots, 0)$ to $(1, \dots, 1)$ given by $d$ linear steps 
in the unit directions $e_1, \dots, e_d$, in this order. The method \texttt{CAxisTensor} computes the 
$k$-th level signature of the canonical axis path in dimension $d$.

The second important case of core tensors is obtained by considering the {\em canonical monomial path} in $\mathbb{R}^d$, 
defined as the path from $(0, \dots, 0)$ to $(1, \dots, 1)$ given by $t\mapsto (t, t^2, \dots, t^d)$. The method  \texttt{CMonTensor} computes the $k$-th level signature of the canonical monomial path in dimension $d$. 

\medskip

In \cite{AFS18}, the authors explore algebraic varieties corresponding to signature tensors of piecewise linear and polynomial paths. To that end, remark that every piecewise linear path in $\RR^d$ with $m$ pieces and every polynomial path of degree $m$ in $\RR^d$ can be represented by a real $d\times m$ matrix, denoted by $X=(x_{i,j})$.
More precisely, if $X_i$ denotes the $i$-th column of $X$, a piecewise linear path can be parametrized by 
\[t\mapsto X_1 + ...+ X_{i-1} + (mt- i + 1) X_i, \, \,\text{ if } \,t\in\left[\frac{i-1}{m},  \frac{i}{m}\right],  \text{ for } i=1,...,m.\]
 On the other hand, each coordinate function $X_i(t)$ of a polynomial path  can be written as  
\[X_i(t) = x_{i,1}t + ...+ x_{i,m}t^m. \]
In both cases, each coordinate $\sigma_{i_1,...,i_k}$ of the tensor $\sigma^{(k)}(X)$  is a homogeneous polynomial of
degree $k$ in the $dm$ unknowns $x_{ij}$. Thus, we can consider them as the homogeneous coordinates of the projective space $\mathbb P^{dm-1}$, and we have that $\sigma^{(k)}$ defines a rational map of degree $k$:
\begin{center}
    \begin{tikzcd}
    &\sigma^{(k)}:\mathbb P^{dm-1} \arrow[r, dashed]  &\mathbb P^{d^k-1} \\
    \end{tikzcd}
\end{center}
We define the \emph{polynomial signature variety} $\mathcal{P}_{d,k,m}$ and the \emph{piecewise linear signature variety} $\mathcal{L}_{d,k,m}$ as the Zariski closure of the corresponding rational maps. 

\bigskip
In this package we provide the necessary tools to create these maps, from which one may compute the dimension, degree and generators of the signature varieties for different values of $d,k,m$.  

\begin{example}
Let us compute the prime ideal of the piecewise linear signature variety $\mathcal{L}_{3,3,2}$.
\begin{verbatim}
    i1 : loadPackage "PathSignatures";
    i2 : d=3; k=3; m=2;
    i3 : R = QQ[a_(1,1)..a_(d,m)];
    i4 : Ma = transpose genericMatrix(R,m,d);
    
                 3      2
    o4 : Matrix R  <-- R
    
    i5 : A2 = wordAlgebra(m, CoefficientRing => R);
    i6 : CAx = CAxisTensor(k, A2)
    
         1   3 1      2 1   2    1   3
    o6 = -Lt  +-Lt Lt  +-Lt  Lt +-Lt
         6  2  2  1  2  2  1   2 6  1
    
    i7 : DAx = Ma * CAx; 
    i8 : parAx = tensorParametrization(DAx); 
    o8 : RingMap R <-- QQ[...] 
    i9 : I = ker parAx;
    i10 : dim I
    o10 = 6
    i11 : degree I
    o11 = 90
    i12 : betti mingens I
    
                 0   1
    o12 = total: 1 163
              0: 1   1
              1: . 162
    
    o12 : BettiTally
\end{verbatim}
These results are consistent with the tenth row of \cite[Table 3]{AFS18}, namely that $\mathcal{L}_{3,3,2}$ is 5-dimensional (projectively), has degree 90 and its ambient (projective) dimension is $3^3 - 2 = 25$ (one dimension is lost because there is one linear generator). The rest of the generators consist of 162 quadrics. 
\end{example}
When symbolic computations are too costly, we can turn to techniques from numerical algebraic geometry. The parametrization map created via \texttt{PathSignatures} can easily be used as input in the \textsc{Macaulay2} package \texttt{NumericalImplicitization}, in particular the methods \texttt{numericalImageDim} and \texttt{numericalImageDegree}.
\begin{example}
We compute numerically the dimension and degree of the polynomial signature variety $\mathcal{P}_{2,4,3}$ and the piecewise linear signature variety $\mathcal{L}_{2,4,3}$ using numerical methods:
    \begin{verbatim}
    i13 : d=2; k=4; m=3;
    i14 : R = CC[a_(1,1)..a_(d,m)]; --for numerical computations   
    i15 : Ma = transpose genericMatrix(R,m,d);
    
                  3      2
    o15 : Matrix R  <-- R
    
    i16 : A2 = wordAlgebra(m, CoefficientRing => R);
    i17 : CMon = CMonTensor(k, A2);
    i18 : DMon = Ma * CMon;
    i19 : parMon = tensorParametrization(DMon,
                                            CoefficientRing => CC);
    o19 : RingMap R <-- CC [...]
    i20 : needsPackage "NumericalImplicitization";
    i21 : numericalImageDim(parMon,ideal 0_R)
    o21 = 6
\end{verbatim}
Taking the same steps, which we omit, for the piecewise linear signature variety we also obtain:
\begin{verbatim}
    i26 : numericalImageDegree(parAx,ideal 0_R,Verbose=>false)
    o26 = 64
\end{verbatim} 
These results agree with the fifth row of \cite[Table 3]{AFS18}: the projective dimension of $\mathcal{P}_{2,4,3}$ and $\mathcal{L}_{2,4,3}$ is 5 and their degrees are 192 and 64, respectively. 
\end{example}

\subsection*{Acknowledgements}
CA and FL acknowledge support from DFG CRC/TRR 388 ``Rough Analysis, Stochastic Dynamics and Related Fields'', Project A04. OR has been supported by European Union’s HORIZON–MSCA-2023-DN-JD programme under the Horizon Europe (HORIZON) Marie Skłodowska-Curie Actions, grant agreement 101120296 (TENORS). We thank the Max Planck Institute for Mathematics in the Sciences, Leipzig for hosting in November 2024 the \textit{Macaulay2 in the Sciences} workshop that started this work, and in particular we thank the organizers Ben Hollering and Mahrud Sayrafi for their support with all technical aspects.

\bibliographystyle{alpha}
\bibliography{bibliography}

\newcommand{\etalchar}[1]{$^{#1}$}
\begin{thebibliography}{AGR{\etalchar{+}}25}

\bibitem[AFS19]{AFS18}
Carlos Am{\'e}ndola, Peter Friz, and Bernd Sturmfels.
\newblock {Varieties of Signature Tensors}.
\newblock {\em Forum of Mathematics, Sigma}, 7:e10, 2019.
\newblock \href{https://doi.org/10.1017/fms.2019.3}{\nolinkurl{doi:10.1017/fms.2019.3}}.

\bibitem[AGR{\etalchar{+}}25]{amendola2025}
Carlos Am\'endola, Francesco Galuppi, \'Angel~David {R\'ios Ortiz}, Pierpaola Santarsiero, and Tim Seynnaeve.
\newblock Decomposing tensor spaces via path signatures.
\newblock {\em Journal of Pure and Applied Algebra}, 229(1):107807, 2025.

\bibitem[BSHW13]{bates2013numerically}
Daniel~J Bates, Andrew~J Sommese, Jonathan~D Hauenstein, and Charles~W Wampler.
\newblock {\em Numerically solving polynomial systems with Bertini}.
\newblock SIAM, 2013.

\bibitem[CGM20]{colmenarejo2020toric}
Laura Colmenarejo, Francesco Galuppi, and Mateusz Micha{\l}ek.
\newblock Toric geometry of path signature varieties.
\newblock {\em Advances in Applied Mathematics}, 121:102102, 2020.

\bibitem[CP20]{colmenarejo2020}
Laura Colmenarejo and Rosa Prei{\ss}.
\newblock Signatures of paths transformed by polynomial maps.
\newblock {\em Beitr. Algebra Geom.}, 61(4):695--717, 2020.

\bibitem[FV10]{Friz_Victoir_2010}
Peter~K. Friz and Nicolas~B. Victoir.
\newblock {\em References}, page 638–651.
\newblock Cambridge Studies in Advanced Mathematics. Cambridge University Press, 2010.

\bibitem[Gal19]{galuppi2019rough}
Francesco Galuppi.
\newblock The rough {V}eronese variety.
\newblock {\em Linear algebra and its applications}, 583:282--299, 2019.

\bibitem[GS02]{grayson2002macaulay2}
Daniel~R Grayson and Michael~E Stillman.
\newblock Macaulay2, a software system for research in algebraic geometry, 2002.

\bibitem[GS24]{galuppi2024rank}
Francesco Galuppi and Pierpaola Santarsiero.
\newblock Rank and symmetries of signature tensors.
\newblock {\em arXiv preprint arXiv:2407.20405}, 2024.

\bibitem[Lyo14]{lyons2014rough}
Terry Lyons.
\newblock Rough paths, signatures and the modelling of functions on streams.
\newblock In {\em Proceedings of the International Congress of Mathematicians}, volume~4, pages 163--184. Kyung Moon Publisher, 2014.

\bibitem[MS21]{michalek2021invitation}
Mateusz Micha{\l}ek and Bernd Sturmfels.
\newblock {\em Invitation to nonlinear algebra}, volume 211.
\newblock American Mathematical Soc., 2021.

\bibitem[Pre23]{preiss2023algebraic}
Rosa Prei{\ss}.
\newblock An algebraic geometry of paths via the iterated-integrals signature.
\newblock {\em arXiv preprint arXiv:2311.17886}, 2023.

\bibitem[PSS19]{pfeffer2019learning}
Max Pfeffer, Anna Seigal, and Bernd Sturmfels.
\newblock Learning paths from signature tensors.
\newblock {\em SIAM Journal on Matrix Analysis and Applications}, 40(2):394--416, 2019.

\bibitem[Ree58]{Ree}
Rimhak Ree.
\newblock Lie elements and an algebra associated with shuffles.
\newblock {\em Annals of Mathematics}, 68(2):210--220, 1958.

\bibitem[Reu93]{reutenauerfree93}
Christophe. Reutenauer.
\newblock {\em Free lie algebras}.
\newblock London Mathematical Society monographs ; new ser., 7. Clarendon Press, Oxford, 1993.

\bibitem[Sul23]{sullivant2023algebraic}
Seth Sullivant.
\newblock {\em Algebraic statistics}, volume 194.
\newblock American Mathematical Society, 2023.

\end{thebibliography}

\end{document}